\documentclass{amsart}

\usepackage{graphicx}

\newtheorem{thm}{Theorem}[section]
\newtheorem{lem}[thm]{Lemma}
\theoremstyle{remark}
\newtheorem{rem}[thm]{Remark}
\theoremstyle{definition}
\newtheorem{dfn}[thm]{Definition}

\begin{document}


\title[Index of rotational constant mean curvature tori in $\mathbb{S}^3$]{Morse index 
of constant mean curvature tori of revolution in the 3-sphere}

\author{Wayne Rossman}
\address{Department of Mathematics, Faculty of Science, Kobe University,
Kobe 657-8501, Japan}
\email{wayne@math.kobe-u.ac.jp}

\author{Nahid Sultana}
\address{Department of Mathematics, Faculty of Science, Kobe University,
Kobe 657-8501, Japan}
\email{nahid@math.kobe-u.ac.jp}

\begin{abstract}
We compute lower bounds for the Morse index and nullity of constant mean 
curvature tori of revolution in the three-dimensional unit sphere.  
In particular, all such tori have index at least five, with 
index growing at least linearly with respect to the number of the 
surfaces' bulges, and the index of such tori can be arbitrarily large.  
\end{abstract}


\maketitle

\section{Introduction}
\label{Introduction}

The {\em Morse index} $\text{Ind}(\mathcal{S})$ of a constant mean 
curvature (CMC) closed 
(compact without boundary) surface $\mathcal{S}$ is a measure of $\mathcal{S}$'s 
degree of instability with respect to area, in 
this sense: We define $\text{Ind}(\mathcal{S})$ 
as the number of negative eigenvalues of $\mathcal{S}$'s 
Jacobi operator $\mathcal L$, where the function space is the $C^\infty$ 
functions from $\mathcal{S}$ to the reals $\mathbb{R}$.  
This is one of two different definitions used in the literature, the 
other being the {\em weak index}, equal 
to the maximal dimension of a vector space of functions given by 
first derivatives of volume-preserving 
variations, all of whose nonzero members come from variations that reduce 
area.  $\text{Ind}(\mathcal{S})$ always equals or is one greater than 
the weak index (\cite{berard-barbosa}, \cite{rossman1}), 
and a CMC surface is {\em stable} exactly when the 
weak index is zero, implying $\text{Ind}(\mathcal{S})$ is then $\leq 1$.  
Various combinations of the two indices are possible: in 
the Euclidean $3$-space $\mathbb{R}^3$, planes have both indices zero, 
spheres are stable with Morse index one, 
and catenoids and Enneper surfaces have both indices 
one \cite{LR} (for noncompact surfaces, definitions of the indices 
must be appropriately adjusted).  Because the two indices differ 
by at most one, we choose to use only $\text{Ind}(\mathcal{S})$ here (like in 
\cite{urbano}), without signifigantly weakening our results.  

The index of minimal surfaces in $\mathbb{R}^3$ has been well 
studied, see \cite{EK}, \cite{FiCo}, \cite{MR}, 
\cite{nayatani} amongst numerous papers.  In the $3$-dimensional 
unit sphere $\mathbb{S}^3$, 
the totally geodesic spheres have index 1 \cite{simons}, the minimal 
Clifford torus has index $5$ and any other closed minimal surface has index $\geq 
6$ (Urbano \cite{urbano}).  

The index of CMC surfaces is generally less accessible -- certainly so in 
the case of $\mathbb{R}^3$ \cite{rossman1}.  As for closed CMC 
surfaces in $\mathbb{S}^3$, it is perhaps known only that the stable CMC surfaces are 
precisely the round spheres \cite{bdcjost}.  
Also, a straightforward computation gives the 
index of flat tori of revolution, which we do here.  
But an analogous result to that of \cite{urbano} is not yet known for 
closed CMC surfaces in $\mathbb{S}^3$.  

Even the Morse index of closed 
CMC surfaces of revolution in $\mathbb{S}^3$ is 
still unknown, so here we find lower bounds for the index of 
such surfaces.  They come in three classes: 
(1) round spheres, (2) tori with two distinct axes of revolution, the 
{\em flat} CMC tori, (3) tori with only one axis of revolution, the 
{\em non-flat} CMC tori.  The index in case (1) is trivial to find, and in 
case (2) is also easily found.  Case (3) is more difficult, and we obtain 
lower bounds for that case.  

\begin{figure}
  \includegraphics[scale=0.3125]{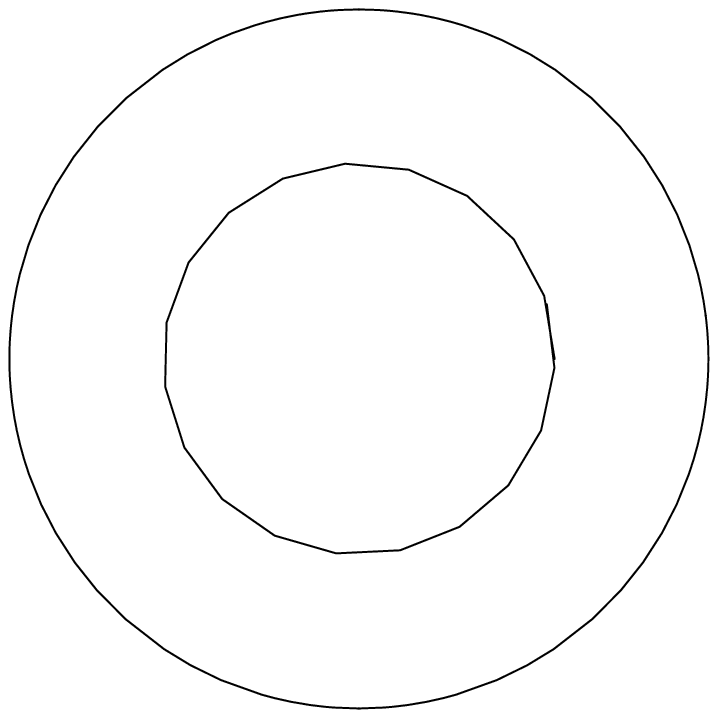}
  \includegraphics[scale=0.424]{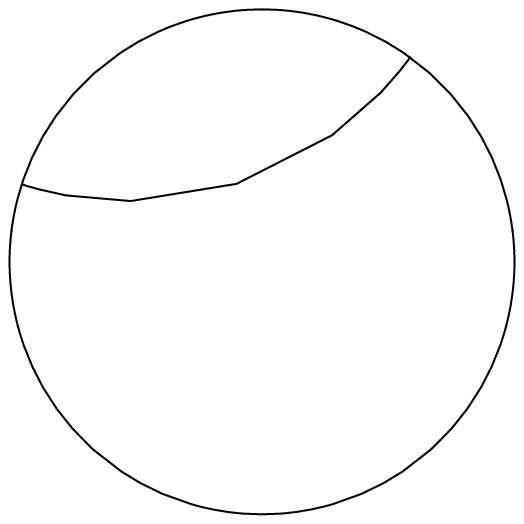}
  \includegraphics[scale=0.325]{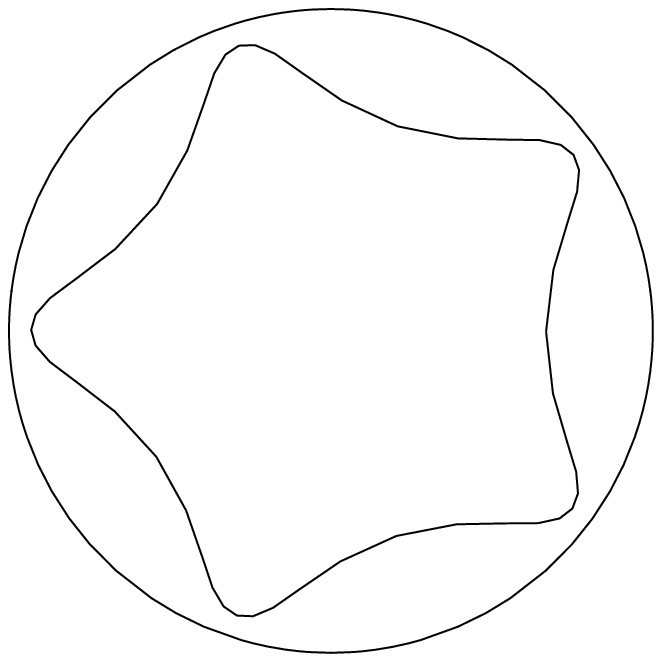}
  \includegraphics[scale=0.344]{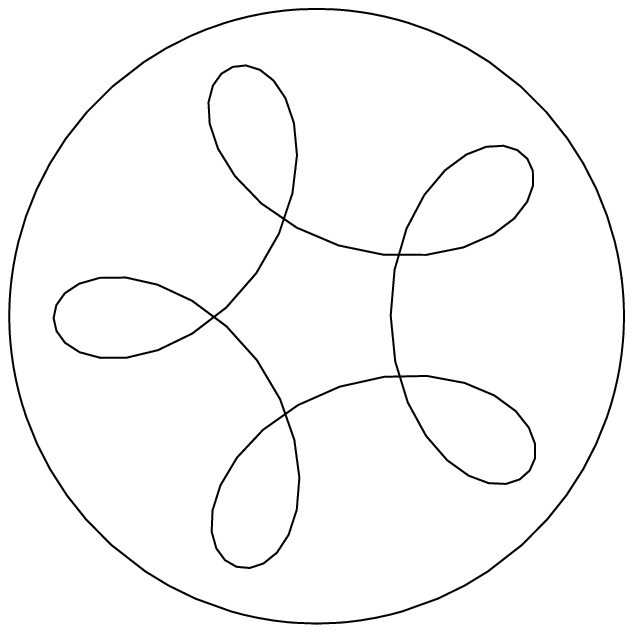}
  \caption{Profile curves of four CMC tori of revolution, in totally geodesic 
           hemispheres having the rotation axis as boundary, are shown.  
           The images are stereographic projections from $\mathbb{S}^3$ to 
           $\mathbb{R}^3 \cup \{ \infty \}$.  The outer circle is the rotation axis, 
           with profile curve inside.  The first curve gives a flat 
           torus, the second a round sphere, the third (resp. fourth) an 
           embedded unduloidal (resp. non-embedded nodoidal) torus with five 
           bulges and five necks.}
  \label{fig1}
\end{figure}

\begin{figure}
  \includegraphics[scale=0.43]{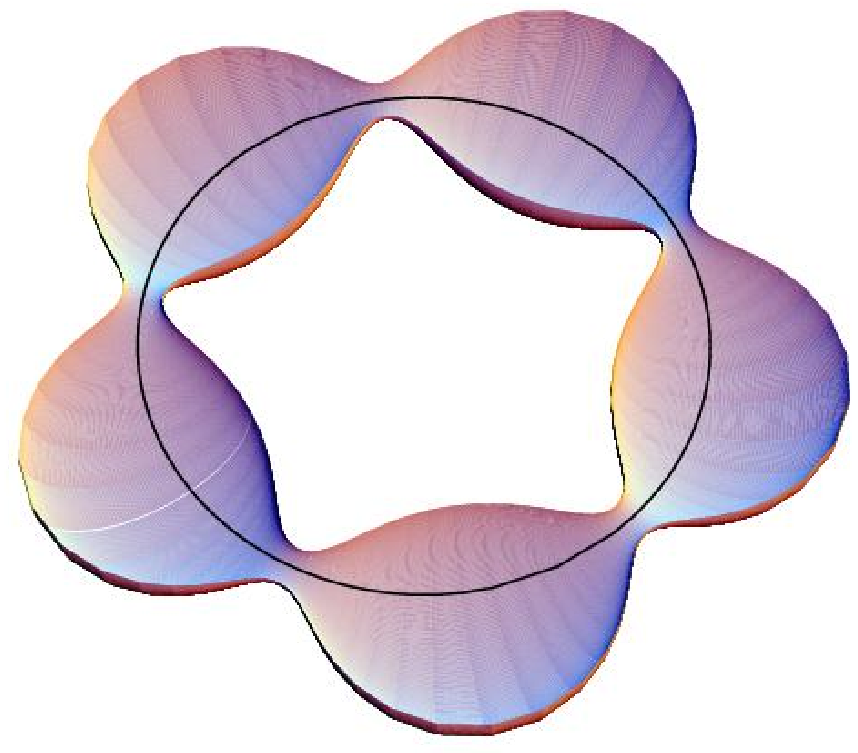}
  \hspace{-1.5in}
  \includegraphics[scale=0.46]{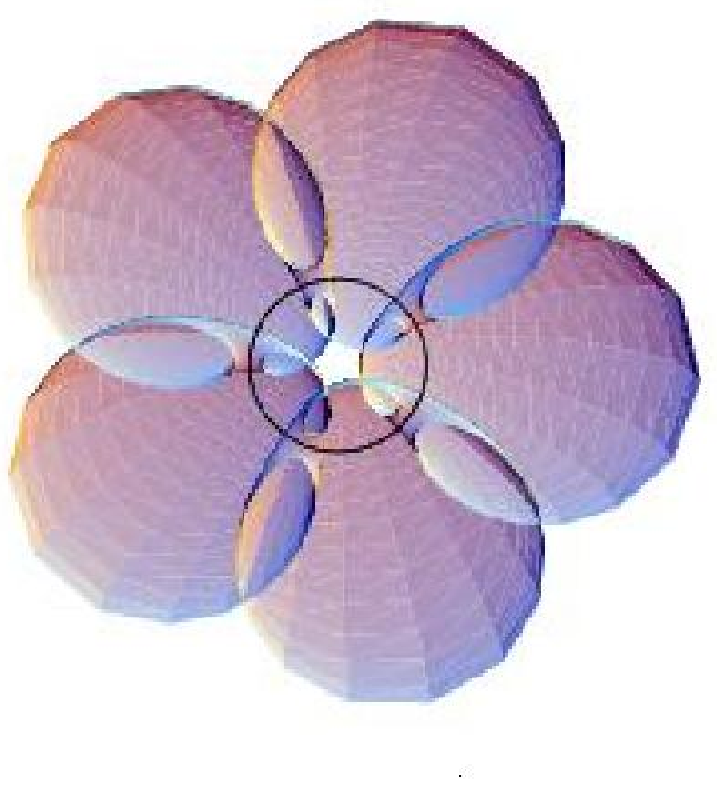}
  \vspace{-0.7in}
  \caption{For the last two profile curves in Figure \ref{fig1}, half of the 
           corresponding surfaces are shown.}
\end{figure}

To state our main result, we give notations: 
For a non-flat CMC torus of revolution $\mathcal{S}$, taking a 
totally geodesic hemisphere whose boundary is the axis of revolution, 
the hemisphere intersects $\mathcal{S}$ along a profile curve with 
an equal finite number of points of maximal and minimal 
distance from the axis \cite{hsiang}, which we call {\em bulges} and 
{\em necks}, respectively.  The nullity is the 
multiplicity of the zero eigenvalue of $\mathcal L$.  
We allow that the mean curvature $H$ can be zero.  

\begin{thm}\label{mainresult} 
Let $\mathcal{S}$ be a closed CMC $H$ surface of revolution in 
$\mathbb{S}^3$.  Then either 
\begin{itemize}
\item $\mathcal{S}$ is a round sphere with Morse index $1$ and nullity $3$, or 
\item $\mathcal{S}$ is a flat torus, and with $b$ equal to 
the greatest integer strictly less than $\sqrt{1+e^{2 \text{arcsinh}(|H|)}}$, the 
index is $3+2b$ and the nullity is either $4$ or $6$, or 
\item $\mathcal{S}$ is a non-flat torus with $k$ bulges and $k$ necks, 
and has index at least $\max(5,2k+1)$ and nullity at least $5$.  
\end{itemize}
\end{thm}

In particular, the index of CMC tori of revolution will always be at least $5$, 
and can indeed be $5$ (when flat and when $H$ is $0$ or close to $0$), and can 
also be arbitrarily large (when flat with $|H|$ large, for example).  

A key ingredient in the proof of the third part of this result is an application of 
Courant's nodal domain theorem, see also \cite{haskins}, 
\cite{rossman1}, \cite{urbano2} for similar applications.  The nodal domain theorem 
can be used to show the index $\geq 2k-1$ by a simpler proof than the 
one here, but the lower bound $\max(5,2k+1)$ is 
much better, in light of \cite{urbano} and 
the facts that the case $k=2$ actually occurs and the case $k=1$ is 
undetermined.  When $H=0$, the nullity is found in \cite{hl}.  Also, the 
first of the three items in Theorem \ref{mainresult} is shown in \cite{simons}.  

Furthermore, when a CMC torus of revolution in $\mathbb{S}^3$ is nodoidal, or is 
unduloidal and wraps at least twice along its axis, the lower bound for the 
index can be improved as in Theorem \ref{improvements} in Section \ref{main-proofs}.  

\section{CMC surfaces of revolution in $\mathbb{S}^3$}
\label{Theparametrizations}

Let $\mathbb{S}^3=\left\{(x_1,x_2,x_3,x_4)\in \mathbb{R}^4 | 
x_1^2 + x_2^2 + x_3^2 + x_4^2=1 \right\}$ inherit its metric $\langle \cdot , 
\cdot \rangle_{\mathbb{S}^3}$ from 
the standard Euclidean $4$-space $\mathbb{R}^4$.  
Up to rigid motions, all CMC surfaces of revolution in 
$\mathbb{S}^3$ can be conformally parametrized (following 
\cite{nahid}, \cite{SKKR}, and secondarily \cite{hsiang}, \cite{mori}, and also 
\cite{hl}, \cite{otsuki}, \cite{yam} when $H=0$) by 
\begin{equation*}
\mathcal{S}_{s,t}(x,y) = 
\left(\text{Re}(\mathcal X)\;,\;\text{Im}(\mathcal Y)\;,\;
\text{Re}(\mathcal Y)\;,\; 
\text{Im}(\mathcal X)\right) : \mathcal{A} \rightarrow \mathbb{S}^3 
\end{equation*}
on the annulus 
$\mathcal{A}=\{ (x,y) \in \mathbb{R}^2 \,| \, (x,y) \equiv (x,y+2 \pi) \}$ (the 
symbol "$\equiv$" denotes an equivalence relation), with
\begin{equation*}\mathcal X=2 e^{i \gamma} \bar{B} \left( 
\frac{c_{-} s_{+} M}{vA\sqrt{CD}}+c_{-}c_{+}\sqrt{\frac{\bar{A}\bar{C}}
{\bar{B}C}}-s_{-}s_{+}\sqrt{\frac{BC}{A\bar{C}}} \right) \;\; , \end{equation*}
\begin{equation*}
\mathcal Y=-\frac{c_{-} c_{+} M}{v\sqrt{ABCD}}-c_{-}s_{+}\sqrt{\frac{\bar{C}}{C}}+
s_{-}c_{+}\sqrt{\frac{C}{\bar{C}}}\;\;,\end{equation*}
where $A=s+te^{2i\gamma}$, $B=se^{2i\gamma}+t$, $C=4ste^{2i\gamma}+v^2$, 
$D=4st+v^2e^{2i\gamma}$, $M=2stv^{\prime}(1-e^{4i\gamma})$, and 
\begin{equation*}
c_{\pm}=\cosh\left(\tfrac{1}{2}(x+iy-g_{\pm})\right)\;\;,\;\;
s_{\pm}=\sinh\left(\tfrac{1}{2}(x+iy-g_{\pm})\right)\;\;,
\end{equation*}
and here 
\begin{equation}\label{closing condition}
(s+t)^2-4 s t \sin^2 \gamma =1/4 \; , \;\; 
g_{\pm}=\int_{0}^{x} 
\frac{2 d\varrho}{1+(4 s t e^{\pm 2 i \gamma})^{-1}v^2(\varrho)} \; , 
\end{equation} 
$s \in \mathbb{R}^+$, $t \in (-s,s] \setminus \{ 0 \}$, 
$\gamma \in (0,\pi/4]$, and $v=v(x)$ solves the ordinary differential equation 
\begin{equation*}
(v^{\prime})^2 = -(v^2-4 s^2)(v^2-4 t^2)\;,\; v(0)=2 t \; , 
\;\;\; v^{\prime} = \tfrac{d}{dx} v \; . 
\end{equation*}
We take $v$ to be the nonconstant periodic solution with values between 
$2|t|$ and $2|s|$ when $|t|<s$, and we take $v$ identically equal to $2 s$ when $s=t$. 

Setting $\rho=16 s^2 t^2 \sin^2 (2 \gamma) v^{-2}$, the metric $ds^2$ and 
Gauss curvature $K$ and mean curvature $H$ for $\mathcal{S}_{s,t}$ are 
\begin{equation*} 
ds^2 = \rho (dx^2+dy^2) \; , \;\; K = -\rho^{-1} (v^2-16 s^2 t^2 v^{-2}) \; , 
\;\; H=\cot(2\gamma) \geq 0 \; . \end{equation*}  

When $t \neq s$, taking $\tau=\sqrt{1-t^2/s^2}>0$, we can write 
$v=2 t/\text{dn}_{\tau} (2sx)$ explicitly via elliptic functions, with 
period (i.e. $v(x)=v(x+x_0)$ $\forall x \in \mathbb{R}$) 
\begin{equation*} 
x_0=\frac{1}{s} \int_0^1 \frac{d\varrho}{\sqrt{(1-\varrho^2)(1- \tau^2 \varrho^2)}} 
\; . \end{equation*} 

By \eqref{closing condition}, we have 
$st \in \left( -(16 \sin^2 \gamma)^{-1},0 \right) \cup 
\left( 0,(16 \cos^2 \gamma)^{-1} \right]$.  When $st$ degenerates to zero, then 
$t=0$ and the profile curve of $\mathcal{S}_{s,t}$ intersects the axis of 
revolution of $\mathcal{S}_{s,t}$, resulting in a round sphere \cite{hsiang}.  
Also by \cite{hsiang}, when $st \neq 0$, we 
know that $\mathcal{S}_{s,t}$ alternates periodically between points of maximum 
distance (bulges) and minimum distance (necks) from the axis of revolution.  In 
analogy to Delaunay surfaces to $\mathbb{R}^3$, we 
have unduloidal surfaces when $st > 0$ and either 
nodoidal or unduloidal surfaces when $st<0$.  We get the minimal Clifford torus 
when $s=t=1/(2 \sqrt{2})$ and $\gamma=\pi/4$, and other flat CMC tori of 
revolution when $s=t$ and $\gamma<\pi/4$.  

Taking $\mu_+=\sqrt{1+16st\sin^2 \gamma}$ and $\mu_- = 
\sqrt{1-16st\cos^2 \gamma}$, and defining 
\begin{equation*} X_+=\frac{-2 \mu_- \sin^2 \gamma + 2 \mu_+ \cos^2 \gamma}{2 + 
(\mu_- + \mu_+) \sin(2\gamma)} \; , \;\;\; 
X_-=\frac{2 \mu_- \sin^2 \gamma + 2 \mu_+ \cos^2 \gamma}{2 - 
(\mu_- - \mu_+) \sin(2\gamma)} \; , \end{equation*} 
and letting $r_+$ and $r_-$ be the distances in $\mathbb{S}^3$ from the axis 
of revolution to the bulges and necks of $\mathcal{S}_{s,t}$, 
we have $r_+ + r_- = 2 \gamma$.  In particular, 
$r_{\pm}=\frac{\pi}{2}-2 \; \text{arctan}(X_{\pm})$.  For 
$st<0$, $\mathcal{S}_{s,t}$ is nodoidal exactly when $r_+ \leq \pi/2$.  

For non-flat $\mathcal{S}_{s,t}$ to close, i.e. to 
become well defined on a torus 
\[ \mathcal T = \{ (x,y) \in \mathbb{R}^2 \,| \, (x,y) \equiv (x,y+2 \pi) \equiv 
(x+x_1,y) \} \] for some $x_1 \in \mathbb{R}^+$, we need 
\begin{equation*} \text{arctan}\left( \frac{2(t+s \cos(2\gamma))}{
\mu_+ \cos^2 \gamma - \mu_- \sin^2 \gamma} \cdot 
 \tan \left(\int_{0}^{x_0} \tfrac{8stv^2(\varrho)\sin(2\gamma) d\varrho}{v^4(\varrho)+
16s^2t^2+8stv^2(\varrho)\cos(2\gamma)} \right) \right) \end{equation*} 
to be a rational multiple of $2 \pi$.  Then $x_1$ will be a 
positive integer multiple of $x_0$.  
In the flat case, $\mathcal{S}_{t,t}$ closes when $x_1=2 \pi/\tan \gamma$.  

\section{The Jacobi operator}
\label{variationofsurfaces}

Endow $\mathcal{T}$ with a 
metric $ds^2=\rho (dx^2+dy^2)$ for $C^\infty$ $\rho=\rho(x,y): \mathcal{T} \to 
\mathbb{R}^+$ (in fact, 
$\rho$ is real-analytic in our application), and let 
$\mathcal{S}:\mathcal{T} \rightarrow \mathbb{S}^3$ be an isometric (hence 
conformal) immersion with mean curvature $H$ and Gauss curvature 
$K$.  When $H$ is constant, $\mathcal{S}$ is critical for a variation problem 
(\cite{berard-barbosa}, \cite{bdcjost}, \cite{blr}, \cite{kp}, 
\cite{silveira}) whose associated {\em Jacobi operator}, or 
{\em stability operator}, is \begin{equation*}
\mathcal L=-\Delta-4-4H^2+2K \; , \end{equation*} 
where $\Delta$ is the Laplace-Beltrami operator of $ds^2$.  The 
function space of $\mathcal L$ is 
taken to be $C^\infty$ functions from $\mathcal{T}$ to 
$\mathbb{R}$.  As $-\mathcal L$ is elliptic and $\mathcal{T}$ is compact, 
it is well known (\cite{berard}, \cite{berger}, \cite{mikhlin}, 
\cite{simons}, \cite{urakawa}) 
that the eigenvalues are real, discrete with finite multiplicities, and 
diverge to $+\infty$.  Since $\Delta$, $K$ and $H$ are all independent of how the 
surface $\mathcal{S}(\mathcal{T})$ is parametrized, clearly the same holds for 
the following $\text{Ind}(\mathcal{S})$ and $\text{Null}(\mathcal{S})$: 

\begin{dfn}\label{defstrongindex}
The {\em index} $\text{Ind}(\mathcal{S})$ of $\mathcal{S}$ 
is the sum of the multiplicities of the negative eigenvalues of 
$\mathcal L$, and the {\em nullity} 
$\text{Null}(\mathcal{S})$ is the multiplicity of the zero eigenvalue of 
$\mathcal L$.  
\end{dfn}

Defining $\hat{\mathcal L}= \rho \mathcal L$, 
the eigenvalues of $\hat{\mathcal L}$ (like for $\mathcal L$) 
form a discrete sequence 
\begin{equation*}
\lambda_1 < \lambda_2 \leq \lambda_3 \leq ...\uparrow +\infty
\end{equation*}
(each considered with multiplicity 1) whose first eigenvalue 
$\lambda_1$ is simple, and the corresponding eigenfunctions 
\begin{equation*}\phi_1,\phi_2,\phi_3,...\in C^{\infty}(\mathcal{T}) \; , \;\;\; 
\hat{\mathcal L} \phi_j=\lambda_j \phi_j \; , \;\; j=1,2,3,... \; , \end{equation*}
can be chosen to form an orthonormal basis with respect to the standard $L^2$ inner 
product $\langle f,g \rangle_{L^2}=\int_{\mathcal T} fg \, dxdy$ for 
functions $f, g : \mathcal{T} \to \mathbb{R}$.  

Furthermore, Courant's nodal domain theorem applies \cite{c}, so 
the number of nodal domains of any eigenfunction associated 
to the eigenvalue $\lambda_j$ is at most $j$.  

$\hat{\mathcal L}$ and $\mathcal L$ will have different 
eigenvalues, but Rayleigh quotient characterizations 
(\cite{berard}, \cite{mikhlin}, \cite{urakawa}) for the eigenvalues 
show that these two operators will give the same index.  Furthermore, as 
the eigenfunctions associated to the zero eigenvalue are the same, 
$\hat{\mathcal L}$ and $\mathcal L$ will also give the same nullity.  
So henceforth we can use either $\mathcal{L}$ or 
$\hat{\mathcal L}$.  In the case that $\mathcal{S}=\mathcal{S}_{s,t}$, then 
\begin{equation*} \hat{\mathcal L} = -\partial_x\partial_x-\partial_y\partial_y-2 v^2 - 
32 s^2 t^2 v^{-2} \; , \end{equation*}  
with function space the $C^\infty$ real-valued functions on $\mathcal T$.  
Later, we also use \begin{equation*} 
\hat{\mathcal L}_0 = -\partial_x\partial_x - 2 v^2 - 
32 s^2 t^2 v^{-2} \; , \end{equation*} now with domain the $C^\infty$ functions 
from the loop $\mathcal T_0 = \{ x \in \mathbb{R} \,| \, x \equiv x+x_1 \}$ to 
$\mathbb{R}$.  The spectrum $\lambda_{1,0} < \lambda_{2,0} \leq 
\lambda_{3,0} \leq ...\uparrow +\infty$ of $\hat{\mathcal L}_0$ has all 
the analogous properties 
as for $\hat{\mathcal L}$, and Courant's nodal domain theorem also applies to 
$\hat{\mathcal L}_0$.  

\section{Proof of the main results}
\label{main-proofs}

For a conformal immersion of a CMC $H$ sphere $\mathcal{S}$ in $\mathbb{S}^3$ 
with Laplace-Beltrami operator $\Delta$, 
the Gauss curvature is $K=1+H^2$ and the Jacobi operator is 
$\mathcal{L}=-\Delta-2(1+H^2)$.  By the canonical 
correspondence \cite{lawson}, the immersion of $\mathcal{S}$ is 
isometric to a conformal immersion of a 
sphere in $\mathbb{R}^3$ of radius $1/\sqrt{1+H^2}$, 
whose first two eigenvalues of $\Delta$ are well known to be 
$0$ and $2(1+H^2)$, with multiplicities $1$ and $3$ respectively (see \cite{berger}, 
for example).  This shows the first item of Theorem \ref{mainresult}.  

When $\mathcal{S}_{t,t}$ is a flat CMC $H$ torus as in Section 
\ref{Theparametrizations}, $\hat{\mathcal L} = - \partial_x \partial_x - 
\partial_y \partial_y - (\cos \gamma )^{-2}$ with eigenvalues 
\begin{equation*}\lambda_{m,n} = m^2 + \alpha n^2 - 1 - 
\alpha \; , \;\;\; \alpha = (\sqrt{1+H^2}-H)^2 \in (0,1] \; , \end{equation*}  
for integers $m,n \geq 0$, with multiplicity $4$, resp. $2$, $1$, 
when $mn>0$, resp. $m+n>0$ and $mn=0$, $m=n=0$.  The region 
$\mathcal{R} = \{ (x,y) \in 
\mathbb{R}^2 \, | \, x^2+\alpha y^2<1+\alpha , x>0, y>0 \}$ is bounded by two 
line segments and part of an ellipse centered about $(0,0)$.  Furthermore, 
$\mathcal R$ contains no points of $\mathbb{R}^2$ with integer coordinates, and 
the points $(1,1)$, $(0,\sqrt{(1+\alpha)/\alpha})$, $(\sqrt{1+\alpha},0)$ and 
$(0,0)$ lie in its boundary.  Now, by summing the multiplicities of the negative 
$\lambda_{m,n}$, we have the index as in the second item of 
Theorem \ref{mainresult}.  The nullity is $6$ when $\sqrt{(1+\alpha)/\alpha}$ 
is an integer, as is $4$ otherwise.  

The third item of Theorem \ref{mainresult} is shown by the following three lemmas: 

\begin{lem}\label{lemmaforproofofmainthm1}
Let $\mathcal{S}(x,y) : \mathcal{A} \to \mathbb{S}^3$ be any conformal 
immersion of revolution with axis $\ell_1 \subset \mathbb{S}^3$, where for 
each $y_0 \in \mathbb{R}$, 
the curve $\hat{c}_{y_0}(x)=\mathcal{S}(x,y_0)$ lies in a unique 
totally geodesic sphere $\mathcal{P}_{y_0}$ of $\mathbb{S}^3$ with 
$\ell_1 \subset \mathcal{P}_{y_0}$.  
By conformality, the angle between any $\mathcal{P}_{y_0}$ 
and $\mathcal{P}_{y_1}$ along $\ell_1$ is $y_0-y_1$, and for each 
$x_0 \in \mathbb{R}$ the curve $\check{c}_{x_0}(y)=\mathcal{S}(x_0,y)$ is 
a circle with center in $\ell_1$.  Let $\mathcal{P}_1 = \mathcal{P}_{\pi/2} 
\supset \hat{c}_{\pm \pi/2}(x)$.  
Let $\vec N$ be a unit normal vector to $\mathcal{S}(x,y)$.  

Let $\mathcal P_2$ be a totally geodesic sphere perpendicular to both 
$\mathcal P_1$ and $\ell_1$, and let $\ell_2$ be a geodesic circle in 
$\mathcal P_2$ intersecting both $\mathcal P_1$ and $\ell_1$ perpendicularly.  
For a Killing field $\mathcal K$ on $\mathbb{S}^3$ produced by constant-speed 
rotation in $\mathbb{S}^3$ about $\ell_2$, $f(x,y) = 
\langle \mathcal K, \vec N \rangle_{\mathbb{S}^3}$ satisfies both 
\begin{itemize}
\item $f(x,y) = u(x) \sin y$ for some function $u(x)$ depending only on $x$, and 
\item if the metric, mean and Gauss curvatures are $\rho (dx^2+dy^2)$, $H$ and 
$K$, and if $\hat{\mathcal L}=-\rho (\Delta + 4 + 4H^2-2K)$, then 
$\hat{\mathcal L} (f(x,y))$ is identically zero.  
\end{itemize}
\end{lem}

\begin{proof}
Let us use the stereographic projection $\mathbb{S}^3=\{ (y_1,y_2,y_3) \, | \, 
y_j \in \mathbb{R} \} \cup \{ \infty \}$, with metric $ds^2 = 4 
(\sum_{j=1}^3 dy_j^2)/(1+\sum_{j=1}^3 y_j^2)^2$.  
(We use this model for $\mathbb{S}^3$ only in this proof.)  
By a rigid motion of $\mathbb{S}^3$, we may assume $\ell_j$ is the 
$y_j$-axis and $\mathcal P_j$ is the $y_jy_3$-plane for $j=1,2$.  
Then $\mathcal{S}$ can be parametrized as 
\begin{equation*} \mathcal{S}(x,y) = 
( \phi(x), \psi(x) \cos y , \pm \psi(x) \sin y ) \end{equation*} for functions 
$\psi(x) \neq 0$ and $\phi(x)$ depending only on $x$.  Then $\vec N$ is of the form 
\begin{equation*} 
\vec N = (a(x), b(x) \cos y , \pm b(x) \sin y ) \end{equation*} for functions 
$a(x)$, $b(x)$ of $x$ with $a^2+b^2=(1+\phi^2+\psi^2)^2/4$.  So 
\begin{equation*} f = u(x) \cdot \sin y \; , \;\;\; u(x) = 
c (b \phi-a \psi)(1+\phi^2+\psi^2)^{-2} \; , \;\;\; c \in \mathbb{R} \setminus 
\{ 0 \} \; , 
\end{equation*} 
and $\hat{\mathcal L} (f)=0$ (see Prop. 2.12 of 
\cite{bdcjost}, for example, or \cite{choe} for when $H=0$).  
\end{proof}

\begin{lem}\label{lemmaforproofofmainthm2}
Let $\mathcal{S}_{s,t}$ be a non-flat CMC torus of revolution 
as in Section \ref{Theparametrizations}, 
with $k$ bulges and associated operator $\hat{\mathcal L}_0$.  
Then $\hat{\mathcal L}_0$ has eigenvalue $-1$ with multiplicity two.  
\end{lem}

\begin{proof}
The axis of $\mathcal{S}_{s,t}$ is $\mathbb{S}^3 \cap \{ x_3=x_4=0 \}$, and 
the point $\mathcal{S}_{s,t}(0,0) \in \mathbb{S}^3 \cap \{ x_2=x_3=0 \}$ 
is a bulge of $\mathcal{S}_{s,t}$.  
Taking two Killing fields of constant-speed rotations of $\mathbb{S}^3$ about 
the geodesic circles $\mathbb{S}^3 \cap \{ x_j=x_3=0 \}$ for $j=1,2$ respectively, 
two corresponding eigenfunctions $u_j(x) \sin y$ ($j=1,2$) of $\hat{\mathcal L}$, 
both with associated eigenvalue zero, are produced as in 
Lemma \ref{lemmaforproofofmainthm1}.  Thus $\hat{\mathcal L}_0 (u_j(x))=-u_j(x)$ 
for $j=1,2$.  Since $\mathcal{S}_{s,t}(0,0)$ is 
a bulge, we have $u_2(0)=0$ and $u_1(0) \neq 0$.  Thus 
$u_1(x)$ and $u_2(x)$ must be independent functions of $x$, and the eigenvalue $-1$ 
of $\hat{\mathcal L}_0$ has multiplicity $\geq 2$.  
Because $\hat{\mathcal L}_0 (u_j)=-u_j$ is a linear second-order ordinary 
differential equation on the domain $\mathcal{T}_0$, $-1$ has multiplicity $\leq 2$, 
and hence exactly $2$.  
\end{proof}

\begin{rem}\label{intermediate-remark}
When $\mathcal{S}_{s,t}$ is nodoidal with $k \geq 2$ 
bulges, it is not difficult to see that 
$u_2$ in the above proof must have at least four zeros.  So, by Courant's nodal 
domain theorem, $\lambda_{j,0}<-1$ for $j=1,2,3$.  
When $\mathcal{S}_{s,t}$ is unduloidal and the projection of a 
profile curve to the axis circle has wrapping number $w \geq 2$ about the circle, 
it follows that $u_2$ must have at least $2 w$ zeros, so 
$\lambda_{j,0}<-1$ for all $j \leq 2w-1$.  
\end{rem}

\begin{lem}\label{lemmaforproofofmainthm3}
Let $\mathcal{S}_{s,t}$ be a non-flat CMC torus of revolution 
with $k$ bulges as in Section 
\ref{Theparametrizations}.  Then $\text{Ind}(\mathcal{S}_{s,t}) \geq \max(5,2k+1)$ and 
$\text{Null}(\mathcal{S}_{s,t}) \geq 5$.  
\end{lem}

\begin{proof}
Recalling the $v(x)$ in Section \ref{Theparametrizations}, a direct computation 
shows that the function 
$u_0(x) = v^\prime/v$ is an eigenfunction of either $\hat{\mathcal L}$ or 
$\hat{\mathcal L}_0$ with 
eigenvalue zero.  (Geometrically, $u_0$ is the oriented length of the 
normal projection of a Killing field produced by constant-speed rotation of 
$\mathbb{S}^3$ about the geodesic of distance $\pi/2$ from the 
axis of $\mathcal{S}_{s,t}$.)  

As $u_0$ is independent of $y$ and has $2k$ nodal domains on 
$\mathcal{T}_0$, Courant's nodal domain theorem implies 
$\hat{\mathcal L}_0$ has at least $2k-1$ negative eigenvalues.  

By Lemma \ref{lemmaforproofofmainthm2}, 
we get two independent eigenfunctions $u_j(x)$ 
of $\hat{\mathcal L}_0$, $j=1,2$, both with eigenvalue $-1$.  
Thus $\lambda_{1,0}<-1$.  Let $u_3(x)$ be an eigenfunction of $\hat{\mathcal L}_0$ 
with eigenvalue $\lambda_{1,0}$.  Hence $u_3$, $u_3 \cos y$, $u_3 \sin y$, $u_1$ and 
$u_2$ are five independent 
eigenfunctions of $\hat{\mathcal L}$ with negative eigenvalues, so 
$\text{Ind}(\mathcal{S}_{s,t}) \geq 5$.  

Furthermore, as $\hat{\mathcal L}_0$ has at least $2k-1$ negative eigenvalues, 
all with associated eigenfunctions independent of $y$, this is true of 
$\hat{\mathcal L}$ as well.  And as 
$u_3 \cos y$ and $u_3 \sin y$ are two more eigenfunctions of 
$\hat{\mathcal L}$ with negative eigenvalues, 
$\text{Ind}(\mathcal{S}_{s,t}) \geq 2k+1$.  

Since $u_0$, $u_1 \cos y$, $u_1 \sin y$, $u_2 \cos y$ and $u_2 \sin y$ are 
five independent eigenfunctions of $\hat{\mathcal L}$ with eigenvalue 
zero, $\text{Null}(\mathcal{S}_{s,t}) \geq 5$.  
\end{proof}

By Remark \ref{intermediate-remark}, when either $\mathcal{S}_{s,t}$ in Lemma 
\ref{lemmaforproofofmainthm3} is nodoidal with $k \geq 2$ bulges or is unduloidal 
with wrapping number $w \geq 2$, we can easily strengthen the proof of 
Lemma \ref{lemmaforproofofmainthm3} to obtain: 

\begin{thm}\label{improvements}
Let $\mathcal{S}$ be a closed CMC $H$ surface of revolution in 
$\mathbb{S}^3$.  
\begin{itemize}
\item If $\mathcal{S}$ is nodoidal with $k \geq 2$ bulges, then 
$\text{Ind}(\mathcal{S}_{s,t}) \geq \max (11,2k+5)$.  
\item If $\mathcal{S}$ is non-flat and 
unduloidal with wrapping number $w \geq 2$ along its axis, 
and with $k$ bulges, then $\text{Ind}(\mathcal{S}_{s,t}) \geq \max (6w-1,2k+4w-3)$.  
\end{itemize}
\end{thm}

\begin{rem}
Using the methods in \cite{rossman2}, one can numerically compute the eigenvalues 
of $\hat{\mathcal L}_0$, and thus the index of, any 
given non-flat torus $\mathcal{S}_{s,t}$, and results of that 
are shown in Table \ref{tablelabel}.  (See \cite{nahid-thesis} for more details.)  
\end{rem}

\begin{figure}
  \includegraphics[scale=0.17]{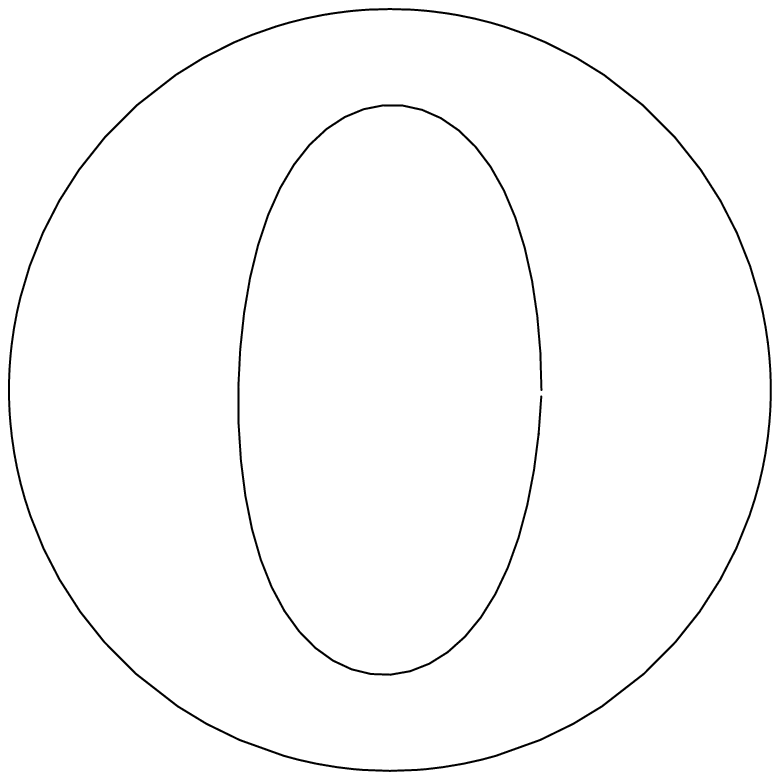}
  \includegraphics[scale=0.17]{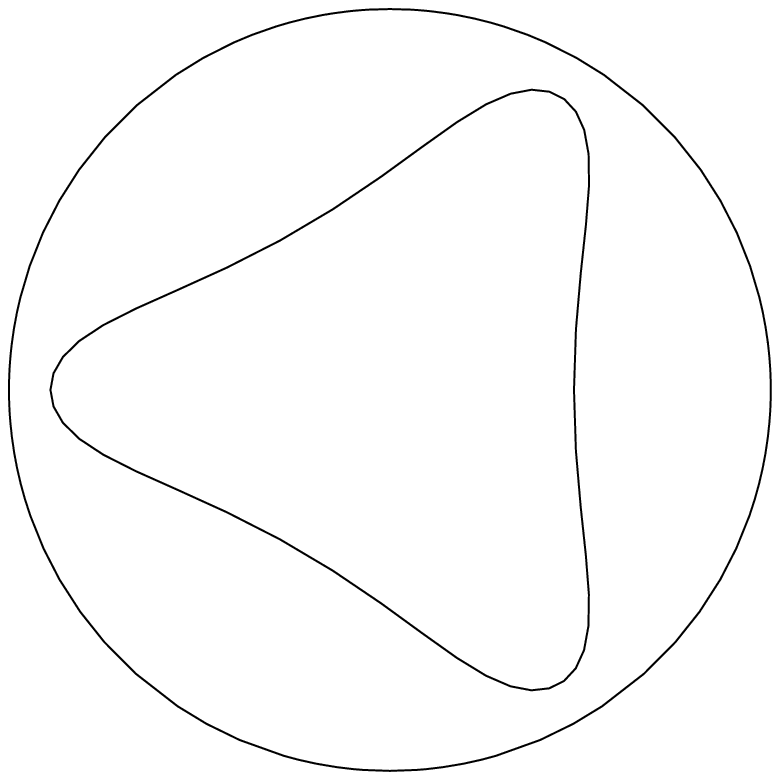}
  \includegraphics[scale=0.17]{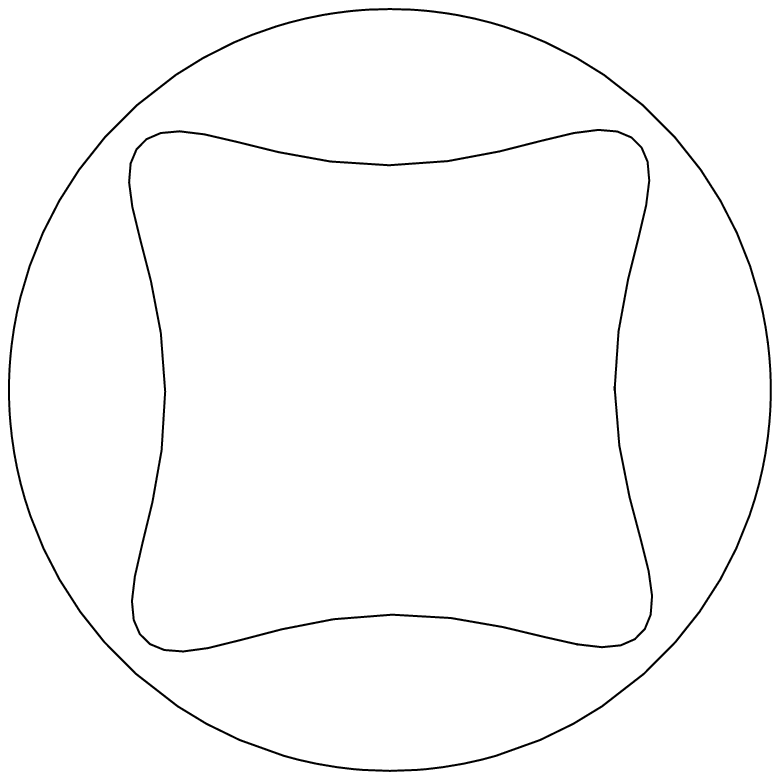}
  \includegraphics[scale=0.2]{5Bulges2.eps}
  \includegraphics[scale=0.17]{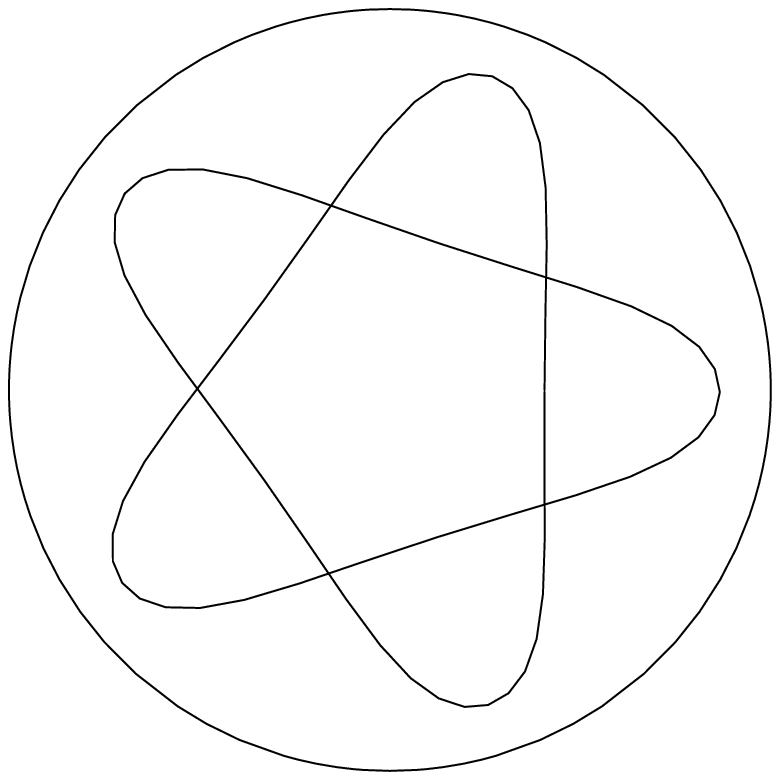}
  \includegraphics[scale=0.17]{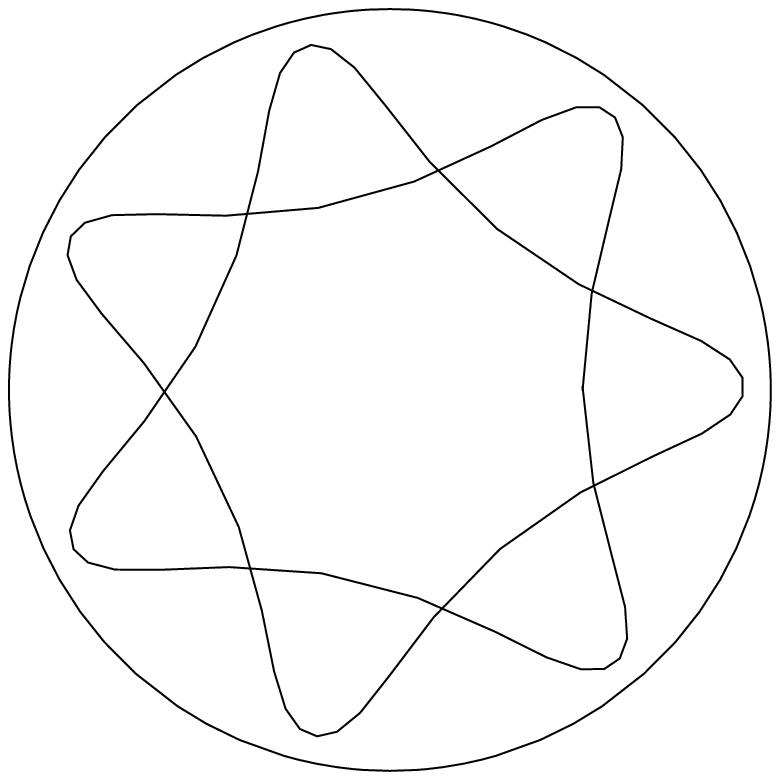}
  \includegraphics[scale=0.17]{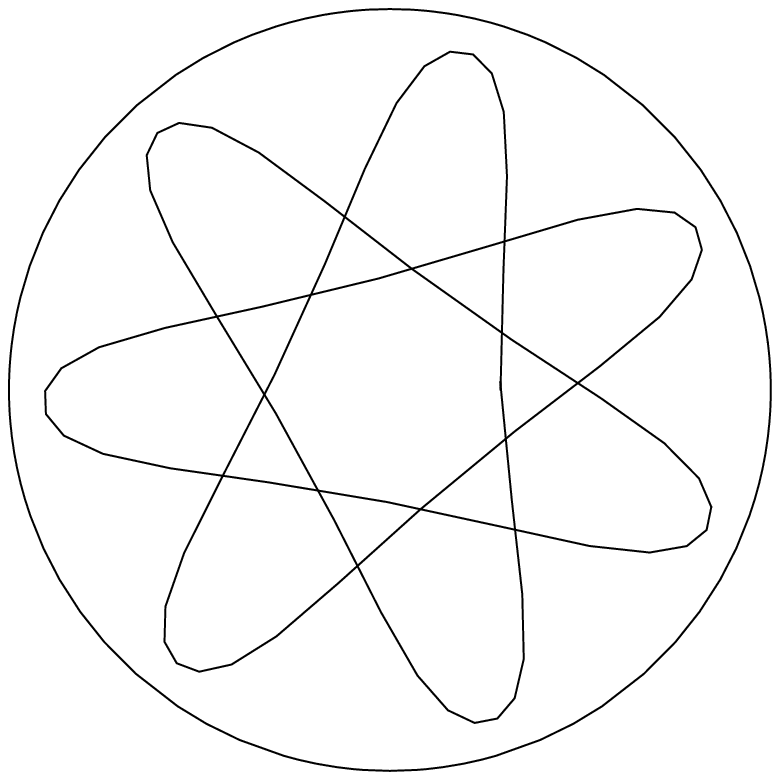}
  \includegraphics[scale=0.17]{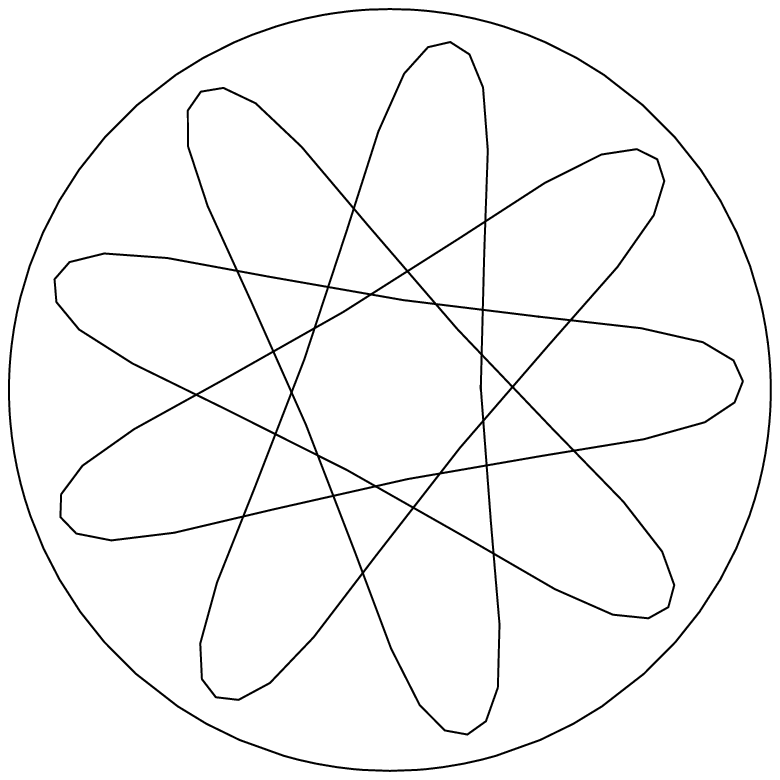}
  \caption{Profile curves for the unduloidal surfaces A, B, C, D, E, F, G and H 
  in Table \ref{tablelabel}.}
  \label{fig2}
\end{figure}

\begin{figure}
  \includegraphics[scale=0.19]{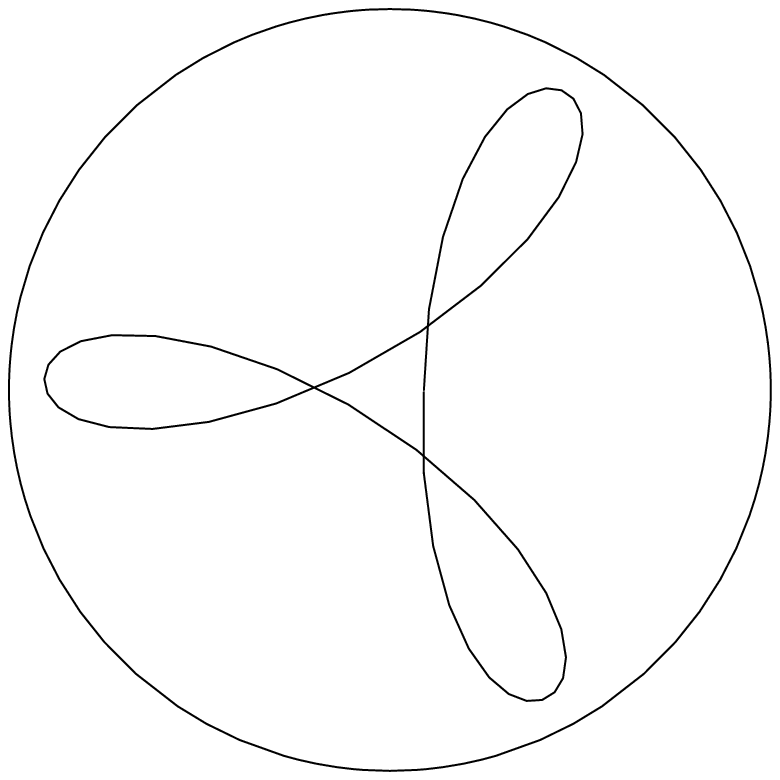}
  \includegraphics[scale=0.276]{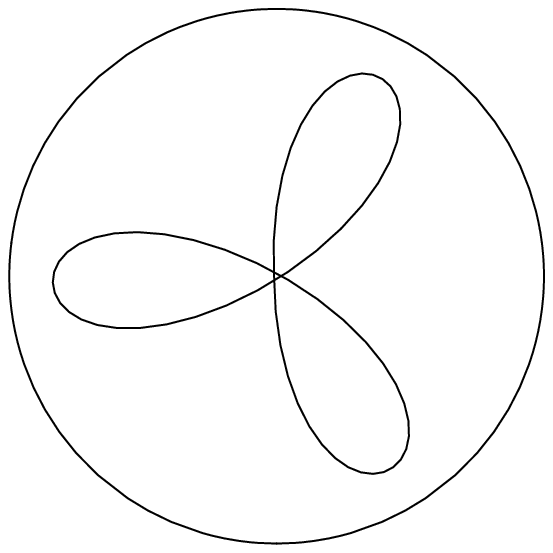}
  \includegraphics[scale=0.19]{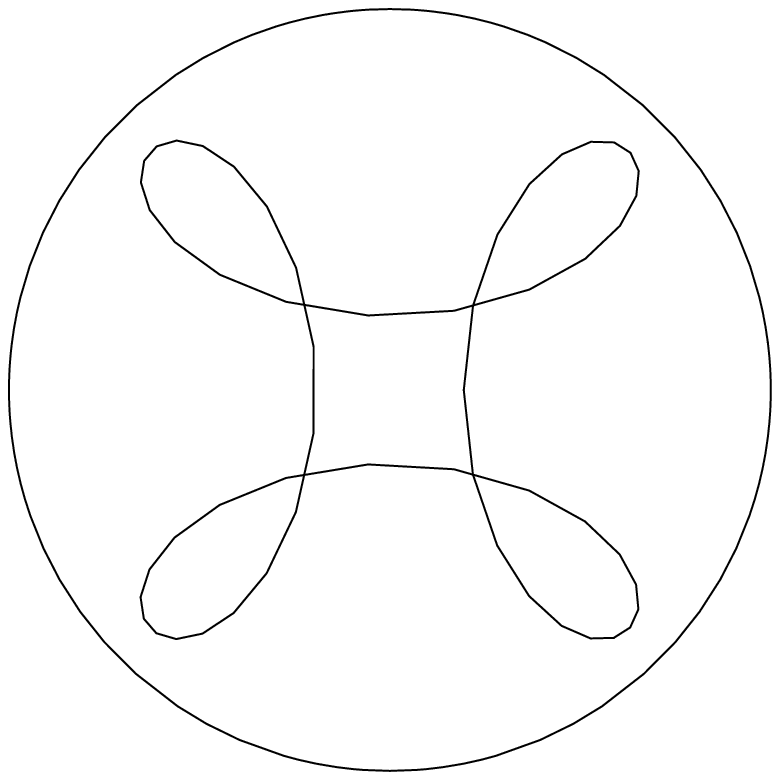}
  \includegraphics[scale=0.19]{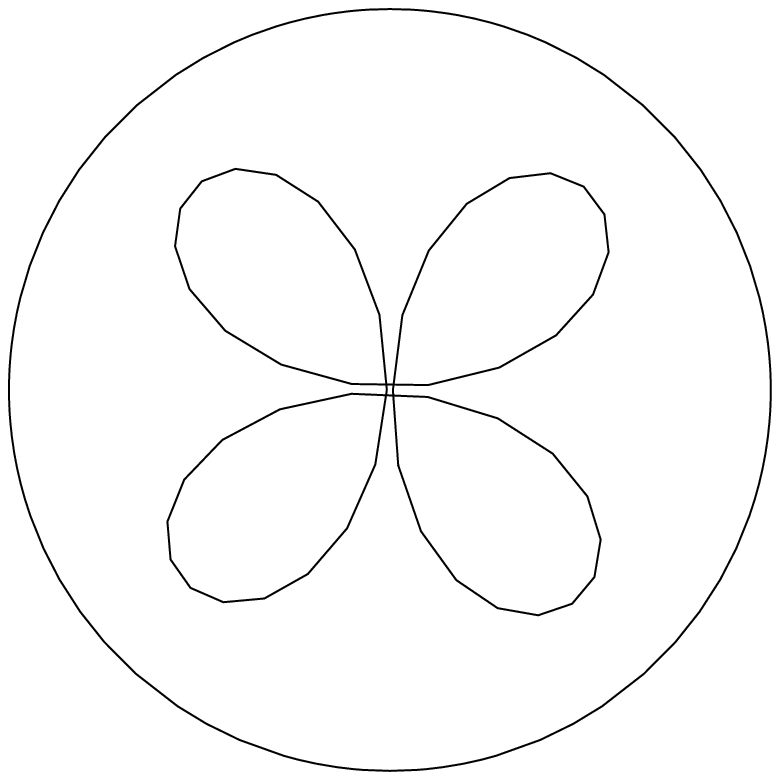}
  \includegraphics[scale=0.237]{5BulgesNod2.eps}
  \includegraphics[scale=0.19]{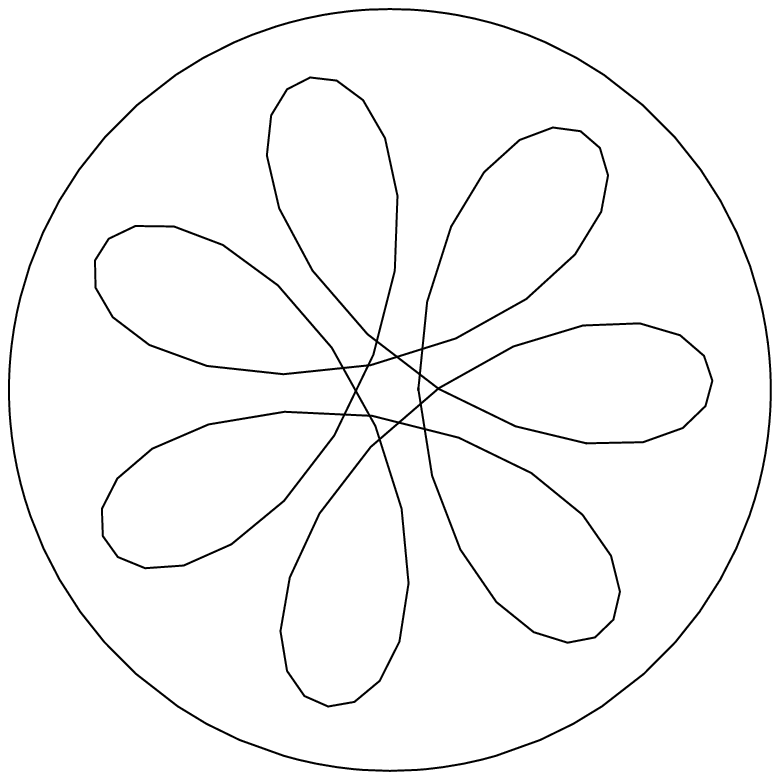}
  \includegraphics[scale=0.19]{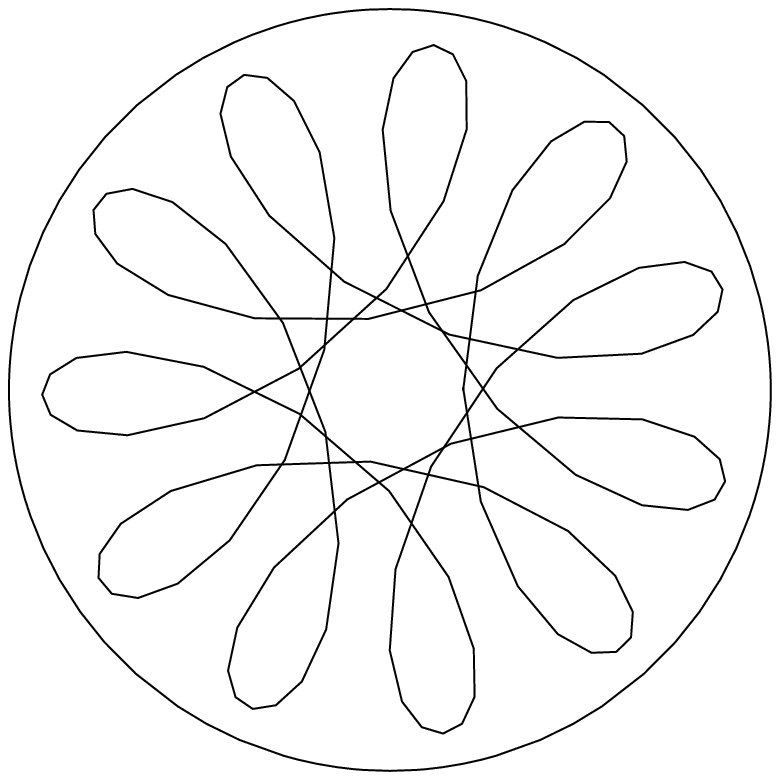}
  \caption{Profile curves for the nodoidal surfaces I, J, K, L, M, N and O 
  in Table \ref{tablelabel}.}
  \label{fig3}
\end{figure}

\begin{rem}
In the proof of Lemma \ref{lemmaforproofofmainthm3}, if 
the eigenvalue of $u_3$ for $\hat{\mathcal L}_0$ 
were strictly less than $-n^2$ for some integer $n \geq 2$, then eigenvalues of 
$\hat{\mathcal L}$ associated to $u_3 \cos (j y)$ and 
$u_3 \sin (j y)$, for integers $j \leq [0,n]$, would be negative, allowing 
us to further strengthen the above results.  However, numerical evaluation shows that 
$\lambda_{1,0} > -2$ for at least the fifteen $\mathcal{S}_{s,t}$ shown in 
Figures \ref{fig2} and \ref{fig3}.  See Table \ref{tablelabel}.  This is related to a 
result in \cite{mori} about the first eigenvalue of the Laplacian on minimal examples. 
\end{rem}

\begin{table}[htb]
{\tiny 
    \begin{center}
    \begin{tabular}{c|c|c|c|c|c|c|c|c|c|c}
      & & & & & & & & numerical & Theorem & Theorem \\
      \text{surf-} & & & & & & & & 
      value & \ref{mainresult}'s lower & \ref{improvements}'s lower \\
     $\text{ace}$ & $s$ & $t$ & $k$ & $w$ & 
$\lambda_{1,0}$ & $\mathcal{B}_-$ & $\mathcal{B}_+$ & of & bound for & bound for \\
  $\mathcal{S}_{s,t}$ & & & & &  & & & $\text{Ind}(\mathcal{S}_{s,t})$ 
  & $\text{Ind}(\mathcal{S}_{s,t})$ & $\text{Ind}(\mathcal{S}_{s,t})$ \\
\hline
      A & $0.4078$ & $0.1583$  & $2$ & $1$& $-1.28$ &$1$ &$1$& $6$& 
           $5$        &$-$ \\ \hline 
      B & $0.4392$ & $0.0811$  & $3$ & $1$&  $-1.08$  &$1$ &$3$& $8$& 
           $7$        &$-$ \\ \hline 
      C & $0.4352$ & $0.0757$  & $4$ & $1$& $-1.04$    &$1$ &$5$& $10$& 
           $9$        &$-$ \\ \hline 
      D & $0.4275$ & $0.0796$  & $5$ & $1$& $-1.02$    & $1$&$7$& $12$& 
           $11$        &$-$ \\ \hline 
      E & $0.4431$ & $0.0881$  & $5$ & $2$& $-1.12$    &$3$ &$5$& $16$& 
            $11$       &$15$ \\ \hline 
      F & $0.4561$ & $0.0559$  & $7$ & $2$& $-1.04$  &$3$ &$9$& $20$& 
            $15$       &$19$ \\ \hline 
      G & $0.4738$ & $0.0527$  & $7$ & $3$ & $-1.11$  &$5$ &$7$& $24$& 
            $15$       &$23$ \\ \hline 
      H & $0.4829$ & $0.0408$  & $9$ & $4$&  $-1.09$ &$7$ &$9$& $32$& 
            $19$       &$31$ \\ \hline 
      I & $0.5112$ & $-0.050$  & $3$ & $1$&  $-1.26$  &$3$ &$1$& $12$& 
            $7$       &$11$ \\ \hline 
      J & $0.5061$ & $-0.089$  & $3$ & $1$&  $-1.40$  &$3$ &$1$& $12$& 
            $7$       &$11$ \\ \hline
      K & $0.5291$ & $-0.068$  & $4$ & $1$& $-1.43$   &$5$ &$1$& $18$& 
            $9$       &$13$ \\ \hline 
      L & $0.5256$ & $-0.155$  & $4$ & $1$&  $-1.85$  &$5$ &$1$& $18$& 
             $9$      & $13$ \\ \hline
      M & $0.5501$ & $-0.095$  & $5$ & $1$& $-1.66$   &$7$ &$1$& $24$& 
            $11$       &$15$ \\ \hline 
      N & $0.5199$ & $-0.087$  & $7$ & $2$&  $-1.47$  &$9$ &$3$& $32$& 
            $15$       &$19$ \\ \hline 
      O & $0.5210$ & $-0.051$  & $11$ & $3$&  $-1.30$  &$15$ &$5$& $52$& 
            $23$       &$27$ 
    \end{tabular}
    \end{center}
}
\caption{Numerical results on the index of the 
CMC non-flat tori of revolution shown in Figures \ref{fig2} and \ref{fig3}.  
There are two $\lambda_{j,0}$ equal to $-1$, and here 
$\mathcal{B}_-$ and $\mathcal{B}_+$ denote the number of $\lambda_{j,0} < -1$ 
and the number of $\lambda_{j,0} \in (-1,0)$ respectively, found 
numerically.  By the numerical 
value of $\text{Ind}(\mathcal{S}_{s,t})$ we mean the value of 
$\text{Ind}(\mathcal{S}_{s,t})=3 \mathcal{B}_- + 
2+ \mathcal{B}_+$ confirmed with numerics, but not yet 
proven with mathematical rigor.}
\label{tablelabel}
\end{table}

\end{document}